\documentclass[a4paper,11pt]{article}
\usepackage{latexsym,amssymb,amsmath,textcomp}
\usepackage[english]{babel}

\newcommand{\E}{\mathbb{E}}

\usepackage{amsthm}
\theoremstyle{plain}
\usepackage{mathtools}

\newtheorem{Theorem}{Theorem}
\newtheorem{Theorem*}{Theorem}
\newtheorem{Remark}{Remark}
\newtheorem{Lemma}{Lemma}

\newtheorem{Proposition}{Proposition}

\usepackage{hyperref}
\usepackage{tikz}

\title{New bounds for perfect $k$-hashing}
\author{Simone Costa\thanks{DICATAM - Sez. Matematica, Universit\`a degli Studi di Brescia, Via
Branze 43, I-25123 Brescia, Italy. email: simone.costa@unibs.it}\ \
and Marco Dalai%
\thanks{DII, Universit\`a degli Studi di Brescia, Via
Branze 38, I-25123 Brescia, Italy. email: marco.dalai@unibs.it}}\begin{document}
\maketitle
\begin{abstract}
Let $C\subseteq \{1,\ldots,k\}^n$ be such that for any $k$ distinct elements of $C$ there exists a coordinate where they all differ simultaneously. Fredman and Koml\'os studied upper and lower bounds on the largest cardinality of such a set $C$, in particular proving that as $n\to\infty$, $|C|\leq \exp(n k!/k^{k-1}+o(n))$. Improvements over this result where first derived by different authors for  $k=4$. More recently, Guruswami and Riazanov showed that the coefficient $k!/k^{k-1}$ is certainly not tight for any $k>3$, although they could only determine explicit improvements for $k=5,6$. For larger $k$, their method gives numerical values modulo a conjecture on the maxima of certain polynomials.

In this paper, we first prove their conjecture, completing the explicit computation of an improvement over the Fredman-Koml\'os bound for any $k$. Then, we develop a different method which gives substantial improvements for $k=5,6$.
\end{abstract}

\section{Introduction}
For positive integers $k\geq 2$ and $n \geq 1$, consider a subset $C\subset \{1,\ldots, k\}^n$ with the property that for any $k$ distinct elements of $C$ there exists a coordinate where they all differ. We call such a set a \emph{perfect $k$-hash code} of length $n$, or simply $k$-hash for brevity. The name is motivated by the idea that if each coordinate of $C$ is interpreted as a $k$-hash function on a set $U$ of cardinality $|C|$, then any $k$ elements of $U$ are hashed onto $\{1,2,\ldots,k\}$ by at least one function. 

Determining the largest possible cardinality of such a set $C$ as a function of $k$ and $n$ is a classic combinatorial problem in theoretical computer science. One standard formulation is to study, for fixed $k$, the grow of the largest possible $|C|$ as $n$ goes to infinity. It is known that $|C|$ grows exponentially in $n$. Then one usually defines the \emph{rate} of the code as \footnote{Here and in the whole paper $\log{x}$ is understood to be in base two.}
\begin{equation}
R=\frac{\log|C|}{n}
\end{equation}
and asks for bounds on the rate of codes of maximal cardinality as $n\to \infty$. This formulation of the problem can also be cast as a problem, in information theory, of determining the zero-error capacity under list decoding for certain channels (REF). 
In this paper we consider upper bounds on $R_k$, defined as the limsup, as $n\to\infty$, of the rate of largest $k$-hash codes of length $n$.

A simple packing argument (see \cite{Elias1}) shows that for all $k\geq 2$ one has $R_k\leq \log(k/(k-1))$. 
For $k=3$, the simplest non-trivial case, this evaluates to $\log(3/2) \approx 0.5850$ and is still the best known upper bound to date (the best lower bound is $1/4\log(9/5)\approx 0.212$).
For $k\geq 4$, the first important result  was derived by Fredman and Koml{\'o}s \cite{FredmanKomlos}, who proved that
\begin{equation}
\label{eq:fredmankomlosRk}
R_k\leq k!/k^{k-1}\,.
\end{equation}
We also refer to \cite{Korner1}, \cite{korner86}, \cite{Korner2} and \cite{Korner3} where the Fredman-Koml{\'o}s bound (and some generalizations to hypergraphs) has been cast using the language of graph entropy and to \cite{nilli} where a simple probabilistic proof has been presented.
Improvements were obtained for $k=4$ in \cite{Arikan1},  \cite{Arikan2} and more recently in \cite{DalaiVenkatJaikumar}, \cite{DalaiVenkatJaikumar2}. The most recent progress we are aware of was obtained in \cite{venkat} where the Fredman-Koml{\'o}s bound is proved to be non-tight for any $k\geq 5$, with an explicit new numerical bound for $k=5,6$. For larger $k$, the authors show that an explicit improvement of the Fredman-Koml{\'o}s bound can be obtained subject to a conjecture on the maxima of certain polynomials. Other recent papers on this topic that deserve to be recalled are \cite{Jaikumar2} where the asymptotic behavior of $R_k$ has been studied and \cite{CostaDalai} where the authors attempt to use the polynomial method to upperbound $R_3$ and they state some limitations of this method.

In this paper we make further progress on this problem. We first prove the conjecture formulated in \cite{venkat} and thus complete their proof of explicit new upper bounds on $R_k$ which beat the Fredman-Koml{\'o}s bound for all $k\geq 5$.
Our main contribution is then to expand on the idea used in \cite{DalaiVenkatJaikumar2} to derive a further improvement for $k=5,6$.

In Section \ref{background} we give a brief summary of the approaches used in \cite{DalaiVenkatJaikumar2} and in \cite{venkat}, upon which we build our contribution. 
In Section \ref{sec:venkat_conj} prove the conjecture stated in \cite{venkat} and give a numerical evaluation of the ensuing bound for $k>6$. In Section \ref{sec:newbounds} we present our improvement for $k=5,6$.

\section{Background}
\label{background}
The bounds presented in \cite{FredmanKomlos}, \cite{Arikan2}, \cite{DalaiVenkatJaikumar2} and \cite{venkat} can all be derived by starting with the following Lemma on graph covering (see \cite{Jkumar}).
\begin{Lemma}[Hansel \cite{hansel}]
Let $K_r$ be a complete graph on $r$ vertices and let $G_1,\ldots,G_m$ be bipartite graphs on those same vertices such that $\cup_{i}G_i=K_r $. Let finally $\tau(G_i)$ represent the number of non-isolated vertices in $G_i$. Then
\begin{equation}
\sum_{i=1}^m \tau(G_i) \geq \log r\,.
\end{equation}
\end{Lemma}
The connection with $k$-hashing comes from the following application. Given a $k$-hash code $C$, fix any $(k-2)$-elements subset $\{x_1,x_2,\ldots, x_{k-2}\}$ in $C$. For any coordinate $i$ let $G_i^{x_1, \ldots, x_{k-2}}$ be the bipartite graph with vertex set $G\setminus\{x_1,x_2,\ldots, x_{k-2}\}$ and edge set
\begin{equation}
\label{eq:FKgraph}
E=\left\{ (v,w) : x_{1,i},x_{2,i},\ldots, x_{k-2,i},v,w \mbox{ are all distinct} \right\}\,.
\end{equation}
Then, since $C$ is a $k$-hash code, we note that $\bigcup_i G_i^{x_1, \ldots, x_{k-2}}$ is the complete graph on $G\setminus\{x_1,x_2,\ldots, x_{k-2}\}$ and so
\begin{equation}
\label{eq:hansel_hash}
\sum_{i=1}^n \tau(G_i^{x_1, \ldots, x_{k-2}})\geq \log(|C|-k+2)\,.
\end{equation}
This inequality can be used to prove upper bounds on $|C|$. Since it holds for any choice of $x_1,x_2,\ldots, x_{k-2}$, one can show that the right hand side is small by proving that left hand side cannot be too large for all possible choices of $x_1,x_2,\ldots, x_{k-2}$. One can either use it for some specific choice or take expectation over any random selection. 

Let $f_{i}$ be probability distribution of the $i$-th coordinate of $C$, that is, $f_{i,a}$ is the fraction of elements of $C$ whose $i$-th coordinate is $a$. Note that the graph in \eqref{eq:FKgraph} is empty if the $x_{1,i},x_{2,i},\ldots, x_{k-2,i}$ are not all distinct. We will say in this case that $x_1,x_2,\ldots,x_{k-2}$ \emph{collide} in coordinate $i$.
Then, we have
\begin{equation}
\tau(G_i^{x_1, \ldots, x_{k-2}})=
\begin{cases}
0 \hspace{2cm} x_1,\ldots,x_k \mbox{ collide in coordinate }i\\
\left(\frac{|C|}{|C|-k+2}\right)\left(1-\sum_{j=1}^{k-2}f_{i,x_{ji}}\right) \hspace{0.5cm}\mbox{otherwise}
\end{cases}
\end{equation}
So, one can make the left hand side in \eqref{eq:hansel_hash} small by either taking a set $x_{1}, \ldots,x_{{k-2}}$ which collide in many coordinates, so forcing the corresponding $\tau$'s to zero, or by taking a set which uses ``popular'' values in many coordinates.

The Fredman-Koml{\'o}s bound is obtained by taking expectation in $\eqref{eq:hansel_hash}$ over a uniform random extraction of $x_1,x_2,\ldots, x_{k-2}$. By linearity of expectation the computation can be performed over each single coordinate. Denoting with $\E$ the expectation, for large $n$ and $|C|$
\begin{align*}
\E [ & \tau(G_i^{x_1, \ldots, x_{k-2}})] \\
&  =\left(1+o(1)\right)\sum_{\stackrel{\text{distinct }}{ a_1,\ldots,a_{k-2}}} f_{i,a_1}f_{i,a_2}\cdots f_{i,a_{k-2}}(1-f_{i,a_1}\cdots-f_{i,a_{k-2}} )
\end{align*}
where the coefficient $o(1)$ is due to sampling without replacement. One can show that the worst-case $f_i$ is the uniform distribution, which gives 
\begin{equation}
\label{eq:fredmankomlostau}
\E [ \tau(G_i^{x_1, \ldots, x_{k-2}})] \leq \frac{k!}{k^{k-1}}\left(1+o(1)\right)\,.
\end{equation}
The procedures used in \cite{DalaiVenkatJaikumar2} and \cite{venkat} are based on the idea that one can also take $x_1,x_2,\ldots, x_{k-2}$ uniformly from a subset  $C'\subset C$ which ensures they collide in all coordinates $i$ in some subset $T \subset \{1,2\ldots, n\}$. Then, if $g_{i,a}$ is the frequency of symbol $a$ in the coordinate $i\notin T$ of $C'$, one has
\begin{align}
\E [ & \tau(G_i^{x_1, \ldots, x_{k-2}})] \nonumber \\
&  =\left(1+o(1)\right)\sum_{\stackrel{\text{distinct }}{ a_1,\ldots,a_{k-2}}} g_{i,a_1}g_{i,a_2}\cdots g_{i,a_{k-2}}(1-f_{i,a_1}\cdots-f_{i,a_{k-2}} )\label{eq:exptaugf}
\end{align}
The worst case $g$ and $f$ here, if taken independently, give in general a value which exceeds the $k!\slash k^{k-1}$ of \eqref{eq:fredmankomlostau}.
In \cite{DalaiVenkatJaikumar2}, for $k=4$, it was shown that one can deal with this by also taking $C'$ randomly from a partition of $C$ (based on the values in positions $i\in T$), thus adding an additional (outer) expectation. In that case $g_i$ is also random and constrained to satisfy $\E[g_i]=f_i$. Using some concavity argument it was shown that under this random selection the bound \eqref{eq:fredmankomlostau} still holds for $i\notin T$, thus gaining on average compared to \cite{FredmanKomlos}. However, for $k>4$ that approach seems infeasible. The idea used in \cite{venkat} is to suppress the random selection of $C'$ and show that one can carefully choose $C'$ so that $x_1, \ldots, x_{k-2}$ collide in a portion of the coordinates large enough to more than compensate the increase in $\E[\tau(G_i^{x_1, \ldots, x_{k-2}})]$ for $i\notin T$ obtained in \eqref{eq:exptaugf} with respect to \eqref{eq:fredmankomlostau}. This leads to a proof that \eqref{eq:fredmankomlosRk} is not tight for all $k>4$. However, explicit numerical improvements were only proved for $k=5,6$, and given for $k>6$ modulo a conjecture  on the optimal value of some polynomials.

In the next two sections we present our contribution. First we prove the conjecture formulated by the authors in \cite{venkat}, thus completing their proof of the new bounds on $R_k$ for all $k$. Then, we prove stronger results for $k=5,6$. Our idea is based on a symmetrization of \eqref{eq:exptaugf} which allow us to resurrect the random selection of $C'$ in an effective way, replacing the concavity argument of \cite{DalaiVenkatJaikumar2} with new bounds on the maxima of some polynomials.

\section{Guruswami-Riazanov bounds}
\label{sec:venkat_conj}

A crucial role in all bounds discussed in this paper is played by the sum appearing in equation \eqref{eq:exptaugf}. We simplify the notation and set, for general probability vectors $g=(g_1,\ldots,g_k)$ and $f=(f_1,\ldots,f_k)$,
\begin{equation}
\psi(g,f) =\sum_{\sigma\in S_k} g_{\sigma(1)}g_{\sigma(2)}\cdots g_{\sigma(k-2)}f_{\sigma(k-1)}\,,
\end{equation}
observing that equation \eqref{eq:exptaugf} can be rewritten as 
\begin{equation}
\label{eq:Expinpsi}
\E [ \tau(G_i^{x_1, \ldots, x_{k-2}})] = \left(1+o(1)\right)\psi(g_i,f_i)
\end{equation}
We can now prove the conjecture stated in \cite{venkat}.

\begin{Proposition}[Conjecture 1 \cite{venkat}]\label{conjecture}
Under the constraints $f_i\geq \gamma,\forall i$, $\psi(g,f)$ attains a maximum in a point $(g,f)$ with vector $f$ of the form $f=(\gamma,\dots,\gamma,1-(k-1)\gamma)$. 
\end{Proposition}
\proof
Since $\psi(g,f)$ is invariant under (identical) permutations on $g$ and $f$, we can study maxima for which $g_k$ is the minimum among the values $g_1,g_2,\dots, g_k$ and show that for those points $f=(\gamma,\dots,\gamma,1-(k-1)\gamma)$. We prove this by considering the components of $f$ one by one. Assume on the contrary that $f_1>\gamma$. Given $\epsilon \leq f_1-\gamma$, set $\tilde{f}=(f_1-\epsilon,f_2\ldots,f_{k-1},f_k+\epsilon)$. Then 
\begin{align*}
\psi(g,\tilde{f}) & = \sum_{\sigma:\sigma(k-1)\neq 1,k} g_{\sigma(1)}g_{\sigma(2)}\dots g_{\sigma(k-2)}f_{\sigma(k-1)}\\ &\qquad +
\sum_{\sigma:\sigma(k-1)=1} g_{\sigma(1)}g_{\sigma(2)}\dots g_{\sigma(k-2)}(f_1-\epsilon)\\
&\qquad +
\sum_{\sigma:\sigma(k-1)=k} g_{\sigma(1)}g_{\sigma(2)}\dots g_{\sigma(k-2)}(f_k+\epsilon)\\
& = \psi(g,f) - \epsilon\cdot \sum_{\sigma:\sigma(k-1)=1} g_{\sigma(1)}g_{\sigma(2)}\dots g_{\sigma(k-2)}\\
& \qquad +
\epsilon \cdot \sum_{\sigma:\sigma(k-1)=k} g_{\sigma(1)}g_{\sigma(2)}\dots g_{\sigma(k-2)}\,.
\end{align*}
Since we assumed $g_1\geq g_k$, 
\begin{equation}
\sum_{\sigma:\sigma(k-1)=1} g_{\sigma(1)}g_{\sigma(2)}\dots g_{\sigma(k-2)}\leq \sum_{\sigma:\sigma(k-1)=k} g_{\sigma(1)}g_{\sigma(2)}\dots g_{\sigma(k-2)}.
\end{equation}
and hence $\psi(g,\tilde{f})\geq \psi(g,f)$. By repeating the above procedure for $f_2$, $f_3,\ldots,f_{k-1}$, we find that indeed $f=(\gamma,\dots,\gamma,1-(k-1)\gamma)$ maximizes $\psi(g,f)$ under the considered constraints whenever $g_k$ is the minimum among $g_1,g_2,\dots, g_k$, and in particular for the optimal $g$ sorted in this way.\endproof

It terms of $g$, it was already shown in \cite{venkat} that assuming the above result one could show that the maximum value of $\psi(g,f)$, under the constraint that $f_i\geq \gamma,\forall i$, is attained at a point $(g,f)$ with $g$ of the form $(\beta,\beta,\ldots,1-(k-1)\beta)$.
Assuming this, it was shown in \cite{venkat} that a new explicit numerical bound can be given on $R_k$ which strictly improves the Fredman-Koml\'os bound for all $k$.
Table \ref{tab:numbounds} gives numerical results\footnote{We believe that due to a minor error in the computation, the bound given for $R_5$ in \cite{venkat} is not really the best possible using their method. We report here the optimal.} for the first values of $k$.

\begin{table}
\centering
\begin{tabular}{|c|c|c|c|c|}
\hline
$k$ & 5 & 6 & 7 & 8 \\
\hline
Bound from \cite{FredmanKomlos} & 0.19200 & 0.092593 & 0.04284 & 0.019227\\
\hline
Bound from \cite{venkat} & 0.19079 & 0.092279 & 0.04279 & 0.019213\\
\hline
\end{tabular}
\caption{Numerical values for the bounds on $R_k$ from \cite{FredmanKomlos} and from \cite{venkat} in light of Proposition \ref{conjecture}. All numbers are rounded upwards.}
\label{tab:numbounds}
\end{table}

\section{Better bounds for small $k$}
\label{sec:newbounds}
In this section we combine insights from both the approaches  of \cite{DalaiVenkatJaikumar2} and \cite{venkat}. Instead of looking at one subcode $C'$, as done in \cite{venkat}, we follow the idea in \cite{DalaiVenkatJaikumar2}. We consider a partition $\{C_{\omega}: \omega \in \Omega\}$ of our $k$-hash code $C$ and randomly select a subcode $C_{\omega}$. Then we randomly extract codewords $x_1,\ldots,x_{k-2}$ from $C_{\omega}$ and bound the expected value in \eqref{eq:exptaugf} over both random code and codewords. At this point, we replace the concavity argument of \cite{DalaiVenkatJaikumar2} with a symmetrization trick combined with new bounds on the maxima of certain polynomials. This procedure leads to the following nontrivial improvement on the rates $R_5$ and $R_6$. 
\begin{Theorem}\label{main}
For $k=5,6$ the following bounds hold 
\begin{itemize}
\item $R_5\leq 0.1697$;
\item $R_6\leq 0.0875$.
\end{itemize}
\end{Theorem}

\subsection{Proof of Theorem \ref{main}}
Here our goal is to find a family of subcodes such that any $k-2$ codewords $x_1,x_2,\dots,x_{k-2}$ of a given subcode $C_{\omega}$ collide in all coordinates of $T=[1,\ell]$ for a carefully chosen value of $\ell$, that is, for any coordinate $t\in T$ there exist $i,j$ such that $x_{i,t}=x_{j,t}$. This will ensure that the coordinates from $T$ contribute $0$ to the LHS of \eqref{eq:hansel_hash}. To do this, we cover all the possible prefixes of length $\ell$; the following lemma can be seen as a special case of the known results on the fractional clique covering number (see \cite{PartialCovering}).
\begin{Lemma}\label{covering}
For any positive $\epsilon$, for $\ell$ large enough, there exists a partition $\Omega$ of $\{1,2,\dots,k\}^{\ell}$ such that:
\begin{enumerate}
\item $|\Omega|\leq \left\lfloor{\left(\frac{k}{k-3}+\epsilon\right)^{\ell}}\right\rfloor$.
\item For all $\omega\in\Omega$ and $i=1,\ldots,\ell$, the $i$-th projection of $\omega$ has cardinality at most $ k-3$. 
\end{enumerate}
In particular, for any $\omega\in\Omega$, any $k-2$ sequences in $\omega$ collide in all coordinates $i=1,\ldots, \ell$.
\end{Lemma}
\proof

For any $i\in [1,k]$, consider the set $A_i=\{i,i+1,\dots,i+(k-4)\}$ , where the sums are performed modulo $k$ in $[1,k]$.
To a string $s=(i_1,\dots,i_{\ell})$ in $[1,k]^{\ell}$ we associate a set $\omega_s=A_{i_1}\times A_{i_2}\times \dots \times A_{i_{\ell}}\subset [1,k]^\ell$. Fix a word $x\in [1,k]^{\ell}$, and choose uniformly at random the string $s$; the probability that $x\not\in \omega_s$ is $1-\left(\frac{k-3}{k}\right)^{\ell}$. Therefore, if we choose randomly $h$ strings $s_1,\dots,s_h$, the probability that $x\not \in (\omega_{s_1}\cup \dots \cup \omega_{s_h})$ is $\left(1-\left(\frac{k-3}{k}\right)^{\ell}\right)^h$. Hence, the expected number of words $x\in [1,k]^{\ell}$ that do not belong to any of the $\omega_{s_1}, \dots, \omega_{s_h}$ is
\begin{align*}
\E(|\{x\in[1,k]^ \ell:\ x\not \in \omega_{s_1}\cup \dots \cup \omega_{s_h}\}|) = k^{\ell}\left(1-\left(\frac{k-3}{k}\right)^{\ell}\right)^h.
\end{align*}
If this value is smaller than $1$, then there exists a choice of $s_1,\dots, s_h$ such that that the family $\{\omega_{s_1},\dots, \omega_{s_h}\}$ covers the whole set $[1,k]^{\ell}$.
This happens whenever
$$k^{\ell}\left(1-\left(\frac{k-3}{k}\right)^{\ell}\right)^h<1 $$
or equivalently
$$h>\frac{-\ell \log{k}}{\log\left(1-\left(\frac{k-3}{k}\right)^{\ell}\right)}\,, $$
which holds for
$$h>\ell \left(\frac{k}{k-3}\right)^{\ell} \frac{\log{k}}{\log{e}}.$$
For $\ell$ large enough, setting $h=\left\lfloor{\left(\frac{k}{k-3}+\epsilon\right)^{\ell}}\right\rfloor$ we have the desired inequality.

Removing possible intersections between the sets $\omega_s$ we obtain a partition of $[1,k]^\ell$ with the desired properties, since condition 2) is satisfied by construction.
\endproof

Let $\Omega=\{\omega_1,\ldots,\omega_h\}$ be a partition of $[1,k]^\ell$ as derived from Lemma \ref{covering} and consider the family of subcodes $C_{\omega_1},\dots, C_{\omega_h}$ of $C$ defined by
$$C_\omega=\{x\in C: (x_1,x_2,\ldots, x_{\ell})\in \omega\}.$$
Clearly, any $k-2$ codewords $x_1,x_2,\dots,x_{k-2}$ of a given subcode $C_{\omega}$ collide in all coordinates of $T=[1,\ell]$. As in \cite{DalaiVenkatJaikumar2}, define a subcode $C_{\omega}$ to be \emph{heavy} if $|C_{\omega}|> n$ and to be \emph{light} otherwise. We can show that, if $\ell$ is not too large, most of the codewords are contained in heavy subcodes.
Indeed, if we consider $\ell$ such that $\left(\frac{k}{k-3}+\epsilon\right)^{\ell}\leq 2^{nR-2\log{n}}$,
that is $\ell\leq\frac{nR-2\log{n}}{\log\left(\frac{k}{k-3}+\epsilon\right)}$, we have that
$$\left|\bigcup_{C_{\omega} \mbox{ is ligth}} C_{\omega}\right|\leq n\left(\frac{k}{k-3}+\epsilon\right)^{\ell}\leq n2^{nR-2\log{n}}= \frac{|C|}{n}.$$
This means that at least a fraction $(1-1/n)$ of the codewords are in heavy subcodes. If we remove from $C$ the light codes, the rate changes by an amount $\frac{1}{n}\log(1-1/n)$, which vanishes as $n$ grows. So, in the following we can assume, without loss of generality, that all the subcodes are heavy. 

We are finally ready to describe our strategy to pick the codewords $x_1,\dots,x_{k-2}$: first we choose a subcode $C_{\omega}$ with probability $\lambda_{\omega}=|C_{\omega}|/|C|$ and then we pick uniformly at random (and without replacement) $x_1,\dots,x_{k-2}$ from $C_{\omega}$. Since those codewords collide in all the coordinates from the set $T=[1,\ell]$, we obtain in \eqref{eq:hansel_hash}:
\begin{align}
\log(|C|-k+2)& \leq\E_{\omega\in \Omega}(\E[\sum_{i\in [\ell+1,n]}\tau(G_i^{x_1,x_2,\dots,x_{k-2}})]) \\
&=\sum_{i\in [\ell+1,n]}\E_{\omega\in \Omega}(\E[\tau(G_i^{x_1,x_2,\dots,x_{k-2}})])\label{eq:sumellton}.
\end{align}
Let again $f_{i}$ be probability distribution of the $i$-th coordinate of $C$, and let $f_{i|\omega}$ be the distribution of the subcode $C_\omega$.
Invoking \eqref{eq:Expinpsi} for the expectation over the random choice of $x_1,\ldots,x_{k-2}$, we can write for $i\in [\ell+1,n]$
\begin{align*}
\E_{\omega\in \Omega} (\E[\tau(G_i^{x_1,x_2,\dots,x_{k-2}})])
=(1+o(1))\sum_{\omega\in \Omega} \lambda_{\omega} \psi(f_{i|\omega},f_i).
\end{align*}
Since $f_i=\sum_{\mu \in \Omega} \lambda_{\mu} f_{i|\mu}$ and $\psi$ is linear in its second variable, we have that
$$
\E_{\omega\in \Omega}(\E[\tau(G_i^{x_1,x_2,\dots,x_{k-2}})]) 
=(1+o(1))\sum_{\omega,\mu\in \Omega}\lambda_{\omega}\lambda_{\mu}\psi(f_{i|\omega},f_{i|\mu})\,.
$$
We exploit now a simple yet effective trick. Since the sum above is symmetric in $\omega$ and $\mu$, we can write
\begin{align}
\E_{\omega\in \Omega} (\E[\tau &(G_i^{x_1,x_2,\dots,x_{k-2}})])  \nonumber\\
&=\left(1+o(1)\right)\frac{1}{2}\sum_{\omega,\mu\in \Omega}\lambda_{\omega}\lambda_{\mu}[\psi(f_{i|\omega},f_{i|\mu})+\psi(f_{i|\mu},f_{i|\omega})].\label{simmetrizzata}
\end{align}
Here, we note that $f_{i|\omega}$ has no relation with $f_{i|\mu}$.
Therefore we can just consider the following polynomial function over two generic probability vectors $p=(p_1,p_2,\dots,p_k)$ and $q=(q_1,q_2,\dots,q_k)$
\begin{align}
\Psi(p;q)& : =\psi(p,q)+\psi(q,p)\nonumber\\
&=\sum_{\sigma\in S_k} p_{\sigma(1)}p_{\sigma(2)}\dots p_{\sigma(k-2)}q_{\sigma(k-1)}+ q_{\sigma(1)}q_{\sigma(2)}\dots q_{\sigma(k-2)}p_{\sigma(k-1)}.\label{eq:defPsi}
\end{align}
Because of \eqref{simmetrizzata}, if $M_k$ is the maximum of $\Psi$ over probabilistic vectors $p$ and $q$, equation \eqref{eq:sumellton} says that
\begin{align*}
\log{|C|} & \leq (1+o(1))\frac{1}{2}(n-\ell)\sum_{\omega,\mu\in \Omega}\lambda_{\omega}\lambda_{\mu}M_k\\
& =(1+o(1))\frac{1}{2}(n-\ell) M_k.
\end{align*}
Recalling that $|C|=2^{n R}$ and taking $\ell=\left\lfloor{\frac{nR-2\log{n}}{\log\left(\frac{k}{k-3}+\epsilon\right)}}\right\rfloor$, we obtain
\begin{align*}
R\leq (1+o(1))\left[1-\frac{R-2\log(n)/n}{\log\left(\frac{k}{k-3}+\epsilon\right)}\right]\frac{M_k}{2}.
\end{align*}
Rearranging the terms, taking $n\to\infty$ first and then $\epsilon\to 0$, we deduce the following proposition.
\begin{Proposition}\label{FromMaxToRate}
Let $M_k$ be the maximum of $\Psi$ over probabilistic vectors $p=(p_1,p_2,\dots,p_k)$ and $q=(q_1,q_2,\dots,q_k)$. Then we have the following upperbound on $R_k$
$$R_k\leq  \left(\frac{2}{M_k}+\frac{1}{\log(k/(k-3))}\right)^{-1}.$$
\end{Proposition}
In the next subsection we will prove that $M_5=\frac{15(48+\sqrt{5})}{1936}\approx 0.389226$ and $M_6=24/125$. This implies Theorem \ref{main}.

\subsection{Bounds on $\Psi$}
The goal of this subsection is to find the maximum of the function $\Psi$ as defined in \eqref{eq:defPsi}. For this purpose we first introduce two lemmas that provide some restrictions on this maximum.
\begin{Lemma}\label{lagrange1}
Let $\bar{p}=(\bar{p}_1,\dots,\bar{p}_k)$ and $\bar{q}=(\bar{q}_1,\dots,\bar{q}_k)$ be two probabilistic vectors. If $(\bar{p};\bar{q})$ is a maximum for $\Psi$ such that $\bar{p}_1,\bar{p}_2,\bar{q}_1,\bar{q}_2$ are nonzero, then also $(\frac{\bar{p}_1+\bar{p}_2}{2},\frac{\bar{p}_1+\bar{p}_2}{2},\bar{p}_3,\dots,\bar{p}_k;\frac{\bar{q}_1+\bar{q}_2}{2},\frac{\bar{q}_1+\bar{q}_2}{2},\bar{q}_3,\dots,\bar{q}_k)$ is a maximum for $\Psi$.
\end{Lemma}
\proof
If $\bar{P}=(\bar{p};\bar{q})$ is a maximum for $\Psi(p;q)$ under the constraints $p_1+p_2+\dots+p_k=1$ and $q_1+q_2+\dots+q_k=1$, then it is a maximum also under the stronger constraints $p_1+p_2=c_1$, $q_1+q_2=c_2$ where $c_1=\bar{p}_1+\bar{p}_2$, $c_2=\bar{q}_1+\bar{q}_2$, and $p_i=\bar{p}_i,q_i=\bar{q}_i$ for $i\in\{3,4,\dots,k\}$. Because of the Lagrange multiplier method this means that:
$$\frac{\partial \Psi}{\partial p_1}\Big|_{\bar{P}}=\frac{\partial \Psi}{\partial p_2}\Big|_{\bar{P}}$$
and
$$\frac{\partial \Psi}{\partial q_1}\Big|_{\bar{P}}=\frac{\partial \Psi}{\partial q_2}\Big|_{\bar{P}}\,.$$
It follows that:
$$(\bar{p}_1-\bar{p}_2)a+(\bar{q}_1-\bar{q}_2)b=0$$
and
$$(\bar{q}_1-\bar{q}_2)d+(\bar{p}_1-\bar{p}_2)c=0$$
where $a=\frac{\partial^2 \Psi}{\partial p_1\partial p_2}\big|_{\bar{P}}$, $b=\frac{\partial^2 \Psi}{\partial p_1\partial q_2}\big|_{\bar{P}}=\frac{\partial^2 \Psi}{\partial q_1\partial p_2}\big|_{\bar{P}}=c$ and $d=\frac{\partial^2 \Psi}{\partial q_1\partial q_2}\big|_{\bar{P}}$.
If we set $\bar{p}_1-\bar{p}_2=x$, $\bar{q}_1-\bar{q}_2=y$, the previous equations became:
$$\begin{cases} ax+by=0;\\
cx+dy=0.
\end{cases}$$
In the case $ad-bc\not=0$ the previous system admits only the solution $x=y=0$ that means $\bar{p}_1=\bar{p}_2$ and $\bar{q}_1=\bar{q}_2$. It is clear that here we have $\bar{p}_1=\frac{\bar{p}_1+\bar{p}_2}{2}=\bar{p}_2$, $\bar{q}_1=\frac{\bar{q}_1+\bar{q}_2}{2}=\bar{q}_2$ and hence the thesis is satisfied. 

Let us assume $ad-bc=0$.
Then there exists a line $L$ of points $P(t)$ such that $P(1)=\bar{P}$, $P(0)=(\frac{\bar{p}_1+\bar{p}_2}{2},\frac{\bar{p}_1+\bar{p}_2}{2},\bar{p}_3,\dots,\bar{p}_k;\frac{\bar{q}_1+\bar{q}_2}{2},\frac{\bar{q}_1+\bar{q}_2}{2},\bar{q}_3,\dots,\bar{q}_k)$ and
$$\frac{\partial \Psi}{\partial p_1}\Big|_{P(t)}-\frac{\partial \Psi}{\partial p_2}\Big|_{P(t)}=\frac{\partial \Psi}{\partial q_1}\Big|_{P(t)}-\frac{\partial \Psi}{\partial q_2}\Big|_{P(t)}=0.$$
It follows that $\Psi(P(t))$ is constantly equal to the value of $\Psi$ in $\bar{P}=P(1)$. Since $(\frac{\bar{p}_1+\bar{p}_2}{2},\frac{\bar{p}_1+\bar{p}_2}{2},\bar{p}_3,\dots,\bar{p}_k,\frac{\bar{q}_1+\bar{q}_2}{2},\frac{\bar{q}_1+\bar{q}_2}{2},\bar{q}_3,\dots,\bar{q}_k)$ belongs to the line $L$, this point is also a maximum for $\Psi$.
\endproof
With essentially the same proof we also obtain the following result.
\begin{Lemma}\label{lagrange2}
Let $\bar{p}=(\bar{p}_1,\dots,\bar{p}_k)$ and $\bar{q}=(\bar{q}_1,\dots,\bar{q}_k)$ be two probabilistic vectors. If $(\bar{p};\bar{q})$ is a maximum for $\Psi$ such that $\bar{p}_1,\bar{p}_2$ are nonzero while $\bar{q}_1=\bar{q}_2=0$ then also $(\frac{\bar{p}_1+\bar{p}_2}{2},\frac{\bar{p}_1+\bar{p}_2}{2},\bar{p}_3,\dots,\bar{p}_k;0,0,\bar{q}_3,\dots,\bar{q}_k)$ is a maximum for $\Psi$.
\end{Lemma}
In the next two lemmas, we will provide some further restrictions on the maximum of $\Psi$ using just some combinatorial arguments.
\begin{Lemma}\label{LemmaAmmazzaCasi1}
We have that:
$$\Psi(0,p_2,\dots,p_k;0,q_2,\dots, q_k)\leq \Psi(0,p_2,\dots,p_k;q_2,0,q_3,\dots, q_k).$$
\end{Lemma}
\proof
Because of the definition, we have that $\Psi(0,p_2,\dots,p_k;0,q_2,\dots, q_k)$ evaluates as 
$$\sum_{\sigma:\ \sigma(k)=1} p_{\sigma(1)}p_{\sigma(2)}\dots p_{\sigma(k-2)}q_{\sigma(k-1)}+ q_{\sigma(1)}q_{\sigma(2)}\dots q_{\sigma(k-2)}p_{\sigma(k-1)}.$$
Similarly, we have that $\Psi(0,p_2,\dots,p_k;q_2,0,q_3,\dots, q_k)$ equals
\begin{align*}
\sum_{\sigma:\ \sigma(k)=1} & p_{\sigma(1)}p_{\sigma(2)}\dots p_{\sigma(k-2)}q_{\sigma(k-1)}+ q_{\sigma(1)}q_{\sigma(2)}\dots q_{\sigma(k)}p_{\sigma(k-1)}+\\
& (k-2)p_2q_2\left(\sum_{\sigma\in Sym(3,\dots,k)} p_{\sigma(3)}\dots p_{\sigma(k-1)}+ q_{\sigma(3)}\dots q_{\sigma(k-1)}\right).
\end{align*}
The claim follows since each term of the last sum is non negative.
\endproof
The following Lemma is in the same spirit of Proposition \ref{conjecture}.
\begin{Lemma}\label{LemmaAmmazzaCasi2}
We have that:
$$\Psi(p_1,\dots,p_{k-3},0,0,0;q_1,q_2,\dots, q_k)\leq \Psi\left(1,0\dots,0;0,\frac{1}{(k-1)},\dots,\frac{1}{(k-1)}\right).$$
\end{Lemma}
\proof
We suppose, without loss of generality that $q_1$ is the minimum among the values $q_1,q_2,\dots,q_{k-3}$. Setting $p=(p_1,\dots,p_{k-3},0,0,0)$ and $q=(q_1,\dots,q_k)$, we have 
\begin{align*}
\Psi(p;q)& =
\sum_{\sigma:\ \sigma(k-1)\not\in \{1,2\}} q_{\sigma(1)}q_{\sigma(2)}\dots q_{\sigma(k-2)}p_{\sigma(k-1)}\\&
+\frac{p_1+p_2}{2}\sum_{\sigma:\ \{\sigma(k-1), \sigma(k)\}=\{1,2\}} q_{\sigma(1)}\dots q_{\sigma(k-2)}\\
&+(p_1q_2+q_1p_2)(k-2)\sum_{\sigma\in Sym(3,\dots,k)} q_{\sigma(3)}\dots q_{\sigma(k-1)}.\end{align*}

Similarly, setting $p'=(p_1+p_2,0,p_3,\dots,p_{k-3},0,0,0)$, we have that:
\begin{align*}
\Psi(p';q)& =
\sum_{\sigma:\ \sigma(k-1)\not\in \{1,2\}} q_{\sigma(1)}q_{\sigma(2)}\dots q_{\sigma(k-2)}p_{\sigma(k-1)}\\&+
\frac{p_1+p_2}{2}\sum_{\sigma:\ \{\sigma(k-1), \sigma(k)\}=\{1,2\}} q_{\sigma(1)}\dots q_{\sigma(k-2)}\\& +(p_1+p_2)q_2(k-2)\sum_{\sigma\in Sym(3,\dots,k)} q_{\sigma(3)}\dots q_{\sigma(k-1)}.
\end{align*}
Since $q_1\leq q_2$ we have that
$$
\Psi(p;q)\leq \Psi(p';q)\,.
$$
Reiterating the previous procedure, since $q_1$ is the minimum among the values $q_1,\dots,q_{k-3}$, we obtain 
\begin{equation}\label{maggiorazione}\Psi(p_1,\dots,p_{k-3},0,0,0;q_1,q_2,\dots, q_k)\leq \Psi(1,0,\dots,0,0;q_1,q_2,\dots, q_k).\end{equation}
Since $q_1$ does not appear in the value of $\Psi(1,0,\dots,0,0;q_1,q_2,\dots, q_k)$, this is certainly maximized for $q_1=0$. Finally, due to the Muirhead's inequality, we obtain that the RHS of \eqref{maggiorazione} is maximized for $q_2=q_3=\dots=q_k=\frac{1}{k-1}$.
\endproof
As a consequence of the previous lemmas, $\Psi$ attains a maximum in a point of one of the following types:
\begin{itemize}
\item[a)] $\left(1,0\dots,0;0,\frac{1}{(k-1)},\dots,\frac{1}{(k-1)}\right)$;
\item[b)] $(1/k,\dots,1/k;1/k,\dots,1/k)$;
\item[c)] $(0,0,\alpha,\dots,\alpha,\beta,\beta;\gamma,\gamma,\delta,\dots,\delta,0,0)$\\ where $(k-4)\alpha+2\beta=1$ and $2\gamma+(k-4)\delta=1$;
\item[d)] $(0,0,\alpha,\dots,\alpha,\beta;\gamma,\gamma,\delta,\dots,\delta,0)$\\
 where $(k-3)\alpha+\beta=1$ and $2\gamma+(k-3)\delta=1$;
\item[e)] $(0,0,1/(k-2),\dots,1/(k-2);\gamma,\gamma,\delta,\dots,\delta)$\\
 where $2\gamma+(k-2)\delta=1$;
\item[f)] $(0,\alpha,\dots,\alpha,\beta;\gamma,\delta,\dots,\delta,0)$\\
 where $(k-2)\alpha+\beta=1$ and $\gamma+(k-2)\delta=1$;
\item[g)] $(0,1/(k-1),\dots,1/(k-1);\gamma,\delta,\dots,\delta)$\\
 where $\gamma+(k-1)\delta=1$.
\end{itemize}
In particular, because of Lemma \ref{LemmaAmmazzaCasi2}, a maximum with three or more $p$-coordinates (resp. $q$-coordinates) equal to zero is also attained in a point of the form $(a)$. Otherwise, there are at most two zero coordinates both for the vector $p$ and for the vector $q$. Due to Lemma \ref{LemmaAmmazzaCasi1}, we can then assume those zeros are in different positions and finally, using Lemma \ref{lagrange1} and \ref{lagrange2}, we obtain the required characterization of the maximum.

For $k=5,6$, we have inspected using Mathematica all cases listed above and determined the maximum explicitly. 
\begin{Theorem}\label{max}
\begin{itemize}
The following hold:
\item for $k=5$, the global maximum of $\Psi$ is $\frac{15(48+\sqrt{5})}{1936}\approx 0.389226$ and is obtained in case $(g)$ with $\delta=1/44(4+\sqrt{5})$ and $\gamma=1-4\delta$;
\item for $k=6$, the global maximum of $\Psi$ is $24/125=0.192$, obtained in case $(a)$.
\end{itemize}
\end{Theorem}
Theorem \ref{main} follows immediately from Theorem \ref{max} and Proposition \ref{FromMaxToRate}.

\begin{Remark}
For $k>6$, the value obtained for $p$ and $q$ as in case $(a)$, which we conjecture to be the true maximum, is too big to improve the known upper bounds on $R_k$. 
\end{Remark}

\section*{Acknowledgements}

This research was partially supported by Italian Ministry of Education under Grant PRIN 2015 D72F16000790001. Helpful discussions with Jaikumar Radhakrishnan and Venkatesan Guruswami are gratefully acknowledged.


\begin{thebibliography}{30}
\bibitem{Arikan1} E. Arikan, 
A Bound on the Zero-Error List Coding Capacity, 
\textit{In Proceedings. IEEE International Symposium on Information Theory} (1993), 152--152.
\bibitem{Arikan2} E. Arikan, 
An upper bound on the zero-error list-coding capacity, 
\textit{IEEE Transactions on Information Theory} \textbf{40} (1994), 1237--1240.
\bibitem{Jaikumar2}
S. Bhandari and J. Radhakrishnan, Bounds on the Zero-Error List-Decoding Capacity of the q/(q-1) Channel, 2018 IEEE International Symposium on Information Theory (ISIT), Vail, CO, 2018, pp. 906-910.
\bibitem{CostaDalai} S. Costa, M. Dalai, A gap in the slice rank of 
$k$-tensors \textit{preprint}, available at arXiv:1905.07355.

\bibitem{DalaiVenkatJaikumar} M. Dalai, V. Guruswami, and J. Radhakrishnan, 
An improved bound on the zero-error listdecoding capacity of the 4/3 channel, 
\textit{IEEE International Symposium on Information Theory (ISIT)} (2017), 1658--1662.
\bibitem{DalaiVenkatJaikumar2} M. Dalai, V. Guruswami, and J. Radhakrishnan, 
An improved bound on the zero-error listdecoding capacity of the 4/3 channel, 
in \textit{IEEE Transactions on Information Theory}, vol. 66, no. 2, pp. 749-756, Feb. 2020
\bibitem{Elias1} P. Elias, 
Zero error capacity under list decoding, 
\textit{IEEE Transactions on Information Theory} \textbf{34} (1988), 1070--1074.
\bibitem{FredmanKomlos} Michael L. Fredman and J{\'a}nos Koml{\'o}s, 
On the Size of Separating Systems and Families of Perfect Hash Functions, 
\textit{SIAM Journal on Algebraic Discrete Methods} \textbf{5} (1984), 61--68.
\bibitem{venkat} V. Guruswami, A. Riazanov, 
Beating Fredman-Komlos for perfect $k$-hashing, 
\textit{Leibniz International Proceedings in Informatics} (2019).
\bibitem{hansel}
G.~Hansel, ``Nombre minimal de contacts de fermature n\'{e}cessaires pour
  r\'{e}aliser une fonction bool\'{e}enne sym\'{e}trique de $n$ variables,''
  \emph{C. R. Acad. Sci. Paris}, pp. 6037--6040, 1964.
\bibitem{Korner1}  J. Korner, 
Coding of an information source having ambiguous alphabet and the entropy of graphs,
\textit{6th Prague Conference on Information Theory} (1973), 411--425.
\bibitem{korner86}
J. Korner, ``Fredman--{K}oml{\'{o}}s bounds and information theory,'' \emph{SIAM
  Journal on Algebraic Discrete Methods}, vol.~7, no.~4, pp. 560--570, 1986.
\bibitem{Korner2} J. Korner and K. Marton, 
New Bounds for Perfect Hashing via Information Theory,
\textit{European Journal of Combinatorics} \textbf{9} (1988), 523--530.
\bibitem{Korner3} J. Korner, Fredman--Koml{\'o}s bounds and information theory,
\textit{SIAM Journal on Algebraic Discrete Methods} \textbf{7} (1986), 560--570.
\bibitem{PartialCovering} R.J. McEliece and E.C. Posner ``Hide and seek, data storage, and entropy,'' \emph{The Annals of Mathematical Statistics}, vol.~42, 5, pp. 1706--1716, 1971.
\bibitem{nilli}
A.~Nilli, ``Perfect hashing and probability,'' \emph{Combinatorics, Probability
  and Computing}, vol.~3, pp. 407--409, 1994.
\bibitem{Jkumar} J. Radhakrishnan, 
Entropy and Counting,
available at: http://www.tcs.tifr.res. in/~jaikumar/Papers/EntropyAndCounting.pdf.
\end{thebibliography}
\end{document}